\documentclass[a4paper,12pt]{article}

\usepackage{amscd,amssymb,amsmath,amsfonts,amsthm, mathrsfs}
\usepackage{stmaryrd}
\usepackage{graphicx}
\usepackage{verbatim,enumerate,paralist}
\usepackage{textcomp}
\usepackage[colorlinks]{hyperref}
\hypersetup{%
  urlcolor=blue,linkcolor=blue,citecolor=blue
}
\usepackage{pdfpages}

\usepackage[T1]{fontenc}
\usepackage[all,2cell,poly,knot,tips]{xy}
\UseAllTwocells
\CompileMatrices

\usepackage[T1]{fontenc}
\usepackage[all,poly,knot,tips]{xy}
\newtheorem{theorem}{Theorem}[section]
\newtheorem{corollary}[theorem]{Corollary}

\newtheorem{proposition}[theorem]{Proposition}
\newtheorem{remark}[theorem]{Remark}

\newtheorem{example}[theorem]{Example}

\newtheorem{question}[theorem]{Question}

\newtheoremstyle{step}{2\bigskipamount}{\medskipamount}{\upshape}{}{\itshape}{. }{ }{\underline{Step~\thestep}}
\theoremstyle{step}

\renewcommand{\thestep}{\arabic{step}}

\newcommand{\lra}{\longrightarrow}

\newcommand{\ra}{\rightarrow}

\newcommand{\Lra}{\Longrightarrow}
\newcommand{\Ra}{\Rightarrow}

\newcommand{\ldual}[1]{\mathord{{\let\nolimits\relax\sideset{^\wedge}{}{#1}}}}
\newcommand{\laction}[2]{\mathord{{\let\nolimits\relax\sideset{^{#1}}{}{#2}}}}
\newcommand{\conj}[2]{\mathord{{\let\nolimits\relax\sideset{^{#1}}{}{#2}}}}

\newcommand{\xra}{\xrightarrow}
\newcommand{\xla}{\xleftarrow}

\def\CA{{\mathscr A}}
\def\CB{{\mathscr B}}
\def\CC{{\mathscr C}}
\def\CD{{\mathscr D}}

\def\CI{{\mathscr I}}

\def\CL{{\mathscr L}}

\def\CP{{\mathscr P}}

\def\CV{{\mathscr V}}

\def\CX{{\mathscr X}}

\def\dd{{\colon}}
\def\ot{{\otimes}}

\DeclareMathAlphabet{\mathbbe}{U}{bbold}{m}{n}

\begin{document}

\author{Branko Nikoli\'c\footnote{This author gratefully acknowledges the support of an International Macquarie University Research Scholarship for the early part of this research.}\ , Ross Street\footnote{This author gratefully acknowledges the support of Australian Research Council Discovery Grants DP160101519 and DP190102432.}\ , and Giacomo Tendas\footnote{This author gratefully acknowledges the support of an International Macquarie University Research Scholarship.}}

\title{Cauchy completeness for DG-categories}
\date{\today}
\maketitle
\noindent {\small{\emph{2010 Mathematics Subject Classification:} 18D20; 18G35; 18E30; 55P42}}
\\
{\small{\emph{Key words and phrases:} chain complex; suspension; mapping cone; differential graded algebra; 
Cauchy completion; DG-category; absolute colimit; total complex.}}
\begin{center}
\small{Centre of Australian Category Theory, Macquarie University 2109, Australia}
\end{center}
\begin{abstract}
\noindent We go back to the roots of enriched category theory and study categories enriched in chain complexes;
that is, we deal with differential graded categories (DG-categories for short).
In particular, we recall weighted colimits and provide examples.
We solve the 50 year old question of how to characterize Cauchy complete DG-categories in terms of existence
of some specific finite absolute colimits.
As well as the interactions between absolute weighted colimits, we also examine the total complex of a chain complex in a DG-category as a non-absolute weighted colimit.  
\end{abstract}

\tableofcontents

\section{Introduction}\label{intro}

The idea of enriching categories so that their homs lie in some category $\CV$ other than $\mathrm{Set}$
goes back to the early days of category theory. After all, additive categories amount to the case where
$\CV$ is the category $\mathrm{Ab}$ of abelian groups. 

The idea of a more general $\CV$ appeared
in the literature \cite{Svarc, LinThesis, LinAuto, KellyTens, MacLBull} in the early 1960s. The detailed
treatment \cite{EilKel1966} by Eilenberg-Kelly gave many examples of $\CV$ including the category of categories 
and the one that initiated their collaboration, namely, the category $\mathrm{DGAb}$ of chain complexes. 
Categories enriched in $\mathrm{DGAb}$ are called DG-categories; 
their theory (limits and colimits in particular) was developed somewhat in \cite{0} including seeing the concepts of
suspension and mapping cones as colimits.  

Since 1969 or earlier, Lawvere was interested in the case where $\CV$ is the reverse-ordered set of non-negative extended real numbers.
Enriched categories are generalized metric spaces.
He was able (see \cite{LawMetric}) to characterize Cauchy completeness for metric spaces in such a way that
the concept applied to categories enriched in very general $\CV$. For ordinary categories, it means that idempotents split.
For additive categories, it means that idempotents split and finite direct sums exist. 

What Cauchy completeness means for DG-categories was therefore an obvious question.
Various weighted (co)limits which exist in Cauchy complete DG-categories were distinguished in \cite{25},
particularly ``mapping cones''.
These were key to the construction for generating the interlocking diagrams which provided complete
invariants for finite diagrams of chain complexes. 
At the end of Section 1 of \cite{25}, a parenthetic statement (not used elsewhere in the paper) 
made a false claim, without proof, about the Cauchy completion of a DG-category.

In the present paper we characterize Cauchy complete DG-categories in terms of existence
of some specific finite absolute colimits.
The case of graded categories (G-categories) is easier and fairly much as expected; we deal with that first. 
We introduce a slightly new kind of absolute DG-colimit, namely, cokernels of protosplit chain maps. 

We conclude with a construction of a new DG-category $\mathrm{DG}\CA$ for each DG-category $\CA$ which restricts to 
the DG-category of chain complexes when $\CA$ is a mere additive (= $\mathrm{Ab}$-enriched) category.
Some properties are discussed.

\section{Chain complexes}\label{chcmpl}

A good early reference for this section is Chapter IV, Section 6 of Eilenberg-Kelly \cite{EilKel1966}.

Let $\CA$ be any additive category (for us this means an $\mathrm{Ab}$-category). 
Let $\mathrm{DG}\CA$ denote the additive category of chain complexes in $\CA$.
An object $A$, called a {\em chain complex} in $\CA$, is a diagram 
\begin{eqnarray}\label{chaincomplex}
\dots \xra{d_{n+2}}A_{n+1}\xra{d_{n+1}}A_{n}\xra{d_{n}}A_{n-1}\xra{d_{n-1}}\dots
\end{eqnarray}
in $\CA$ satisfying $d_n\circ d_{n+1} = 0$.
We sometimes write $d$ for $d_n$ and $d^A$ for $d$ when there is chance of ambiguity; these are called the {\em differentials} of $A$. 
A morphism $f\dd A\to B$, called a {\em chain map}, is a morphism of diagrams; that is,
a family $$f=(f_n\dd A_n\to B_n)_{n\in \mathbb{Z}}$$ of morphisms in $\CA$
satisfying $d^B\circ f_{n+1}=f_n \circ d^A$ for all $n\in \mathbb{Z}$. 

Let $\mathrm{G}\CA$ denote the additive category of $\mathbb{Z}$-graded objects in $\CA$; we regard it as
the full additive subcategory of $\mathrm{DG}\CA$ consisting of chain complexes with all differentials equal to zero.

The {\em suspension}\footnote{In some literature, $\mathrm{S}^{n}A$ is denoted by $A[n]$.} 
$\mathrm{S}A$ 
of a chain complex $A$ in $\CA$ is defined by $(\mathrm{S}A)_n = A_{n-1}$
with $d^{\mathrm{S}A}=-d^A$. On the chain map $f\dd A\to B$ it is defined by $(\mathrm{S}f)_n = f_{n-1}$.
It is sometimes convenient to regard the differentials as components of a chain map $d^A : A\to \mathrm{S}A$.
Suspension restricts to $\mathrm{G}\CA$.

A chain complex in the opposite category $\CA^{\mathrm{op}}$, called a {\em cochain complex} in $\CA$, is traditionally 
denoted by a diagram 
\begin{eqnarray}\label{superscr}
\dots \xla{d^{n+2}}A^{n+1}\xla{d^{n+1}}A^{n}\xla{d^{n}}A^{n-1}\xla{d^{n-1}}\dots
\end{eqnarray}
in $\CA$ using superscripts rather than subscripts. The virtue of indexing chain complexes by $\mathbb{Z}$
rather than $\mathbb{N}$ is that the cochain complex $A = A^{-}$ is traditionally identified with the chain complex
$A = A_{-}$ defined by 
\begin{eqnarray}\label{cochain}
A_n = A^{-n} \ \text{ and } \ (A_n\xra{d_n}A_{n-1}) = (A^{-n}\xra{d^{-n+1}}A^{-n+1}) \ .
\end{eqnarray}
In fact this defines an additive isomorphism $\mathrm{DG}\CA \cong \mathrm{DG}(\CA^{\mathrm{op}})^{\mathrm{op}}$. 

When $\CA=\mathrm{Ab}$, we have the additive category $\mathrm{DGAb}$. 
The monoidal structure on $\mathrm{DGAb}$ of interest has tensor product 
denoted $A\ot B$ and is defined in terms of the usual tensor product of abelian groups as 
\begin{eqnarray}
& (A\ot B)_n=\sum_{p+q=n}{A_p\ot B_q} \nonumber \\
& d(a\ot b)=da\ot b+(-1)^pa\ot db, \  a\in A_p \ .
\end{eqnarray}
We identify abelian groups with chain complexes $A$ in which $A_n=0$ for all $n\ne 0$;
using that convention, the unit for chain complex tensor product is $\mathbb{Z}$.
Indeed, the monoidal category $\mathrm{DGAb}$ is closed with internal hom $[B,C]$
defined by 
 \begin{eqnarray}\label{chainhom}
& [B,C]_n=\prod_{r-q=n}{\mathrm{Ab}(B_q, C_r)} \nonumber \\
& (df)_qb=d(f_qb)-(-1)^nf_{q-1}db, \  f\in [B,C]_n, \ b\in B_q \ .
\end{eqnarray}

As mentioned in \cite{EilKel1966}, there is a unique symmetry $\sigma = \sigma_{A,B} \dd A\ot B\to B\ot A$
defined by
\begin{eqnarray}
\sigma(a\ot b)=(-1)^{pq} b\ot a,  \  a\in A_p, \  b\in B_q \ .
\end{eqnarray}

We can regard the category $\mathrm{GAb}$ of graded abelian groups as the full subcategory
of $\mathrm{DGAb}$ consisting of the chain complexes $A$ with all $d=0$. The inclusion
functor has right adjoint 
\begin{eqnarray}\label{Zed}
\mathrm{Z}\dd \mathrm{DGAb} \to \mathrm{GAb}
\end{eqnarray}
defined by
$(\mathrm{Z}A)_n= \mathrm{Z}_nA = \{a\in A_n \dd da=0\}$, the group of {\em $n$-cycles} of $A$.
Moreover, the inclusion
has left adjoint $\mathrm{Z}'\dd \mathrm{DGAb} \to \mathrm{GAb}$ defined by
$(\mathrm{Z}'A)_n= \mathrm{Z}'_nA = A_n/\mathrm{im}d^A$.
The image of the canonical natural transformation $\mathrm{Z}\Rightarrow \mathrm{Z}'$
is the {\em homology functor} $\mathrm{H}\dd \mathrm{DGAb} \to \mathrm{GAb}$.
So $(\mathrm{H}A)_n = \mathrm{H}_nA = \mathrm{Z}_nA/\mathrm{im}d^A$.

Notice that $\mathrm{GAb}$ is closed under the tensor product, unit object and internal hom of $\mathrm{DGAb}$,
and so is a symmetric closed monoidal category in such a way that the inclusion is symmetric, strong monoidal,
and strong closed.
It follows that $\mathrm{Z}$ is monoidal and $\mathrm{Z}'$ is opmonoidal.
From the diagram
\begin{equation*}
\xymatrix{
\mathrm{im}d^A\otimes \mathrm{im}d^B\ar[r]^-{} &\mathrm{Z}A\otimes \mathrm{Z}B \ar[r]^-{} \ar[d]_-{\phi} & \mathrm{H}A\ot \mathrm{H}B \ar[d]^-{\phi}\ar[r]^-{} & 0 \\
& \mathrm{Z}(A\ot B) \ar[r]_-{} &  \mathrm{H}(A\ot B) \ar[r]_-{} & 0}
\end{equation*}
with exact rows, we see that $\mathrm{H}$ has monoidal structure induced by that of $\mathrm{Z}$.   

The conservative functor $\mathrm{U}\dd \mathrm{DGAb}\to \mathrm{GAb}$ which forgets all the differentials $d$ (replaces
them by $0$ if you prefer) is also symmetric, strong monoidal (``tensor preserving''), and strong closed (``internal hom preserving''). 
It also has both adjoints.
The left adjoint $\mathrm{L}\dashv \mathrm{U}$ is defined by $(\mathrm{L}X)_n=X_{n+1}\oplus X_n$
with 
\begin{eqnarray}
d= \scriptsize{\left[ 
\begin{array}{ccc}
0 & 1 \\
0 & 0 
\end{array} \right]} \ .
\end{eqnarray} 
The right adjoint $\mathrm{R}\vdash \mathrm{U}$ is defined by $(\mathrm{R}X)_n=X_n\oplus X_{n-1}$
with $d$ given by the same $2\times 2$ matrix.
These facts have a number of pleasant consequences,
the most obvious is that $\mathrm{U}$ preserves and reflects all limits and colimits. 

As an easy application of Beck's monadicity theorem \cite{CWM}, the functor $\mathrm{U}$ is
both monadic and comonadic. The monad $\mathrm{T}= \mathrm{U}\mathrm{L}$ is opmonoidal. 
The comonad $\mathrm{G}= \mathrm{U}\mathrm{R}$ is monoidal. To see a chain complex
$A$ as an Eilenberg-Moore $\mathrm{T}$-algebra, we take the action 
$\alpha \dd \mathrm{T}\mathrm{U}A \to \mathrm{U}A$ to have components
$[d,1]\dd A_{n+1}\oplus A_n\to A_n$. 
There is a $\mathrm{U}$-split coequalizer
\begin{eqnarray}\label{Usplitcoeq}
\xymatrix{
& \mathrm{LULU} \ar@<+1.2ex>[rr]^{\varepsilon \mathrm{LU}}  \ar@<-1.5 ex>[rr]^{\mathrm{LU} \varepsilon} &&   \mathrm{LU}  \ar@{->}[rr]^{\varepsilon}  && 1_{\mathrm{DGAb}} 
}
\end{eqnarray}
which evaluates at $A$ to give chain maps with components
\begin{eqnarray}\label{evsplitcoeq}
\xymatrix{
& A_{n+2}\oplus A_{n+1}\oplus A_{n+1}\oplus A_{n} \ar@<+1ex>[rr]^{\phantom{AAAAA}\beta}  \ar@<-1.5 ex>[rr]^{\phantom{AAAAA}\gamma} &&   A_{n+1}\oplus A_{n}  \ar@{->}[rr]^{\phantom{AA}\alpha}  && A_n 
}
\end{eqnarray}
where 
$\alpha =
\scriptsize{\begin{bmatrix}
d & 1 
\end{bmatrix}}$,
$\beta =
\scriptsize{\begin{bmatrix}
0 & 1 & 1 & 0 \\
0 & 0 & 0 & 1 
\end{bmatrix}}$, 
$\gamma =
\scriptsize{\begin{bmatrix}
d & 1 & 0 & 0 \\
0 & 0 & d & 1 
\end{bmatrix}} $ .

\begin{remark}\label{Hopf}
{\em Indeed, there is a commutative and cocommutative Hopf graded ring $\mathrm{T}\mathbb{Z}$, 
which recaptures the monad $\mathrm{T}$ by tensoring with that ring and determines the monoidal structure on 
$\mathrm{DGAb}$ (as the category of $\mathrm{T}\mathbb{Z}$-modules) via the graded coring 
structure. As this is an aspect of a fairly general setting and since it is not required for our
present purpose, we redirect the reader to \cite{Pareigis81, PMcC2000, 132}.}  
\end{remark}

\section{Constructions on chain complexes}\label{Cocc}

As a category of additive functors into $\mathrm{Ab}$, finite direct sums in $\mathrm{DGAb}$ 
exist and are performed pointwise.  

As mentioned, the suspension $\mathrm{S}A$ 
of a chain complex $A$ is defined by $(\mathrm{S}A)_n = A_{n-1}$
with $d^{\mathrm{S}A}=-d^A$. 
Suspension restricts to $\mathrm{GAb}$ and we have the obvious equality $\mathrm{R} = \mathrm{LS}$;
indeed, we have an infinite string of adjunctions
\begin{eqnarray}\label{stringadjns}
\mathrm{LS}^n\dashv \mathrm{S}^{-n}\mathrm{U}\dashv \mathrm{LS}^{n+1} \ .
\end{eqnarray}

\begin{proposition}\label{Sten}
For chain complexes $A$ and $B$, the equality
\begin{eqnarray*}
\mathrm{S}(A\otimes B) = \mathrm{S}A\otimes B 
\end{eqnarray*}
is a natural isomorphism in $\mathrm{DGAb}$.
Taking $A=\mathbb{Z}$ yields
\begin{eqnarray*}
\mathrm{S}B \cong \mathrm{S}\mathbb{Z}\otimes B \ . 
\end{eqnarray*}
Similarly, for chain complexes $B$ and $C$, the equalities
\begin{eqnarray*}
\mathrm{S}[B,C] = [\mathrm{S}^{-1}B,C] = [B,\mathrm{S}C] 
\end{eqnarray*}
are natural isomorphisms in $\mathrm{DGAb}$.
\end{proposition}
\begin{proof}
Indeed we do have a graded equality 
$$\mathrm{S}(A\otimes B)_n= (A\otimes B)_{n-1}=\sum_{p+q=n-1}{A_p\otimes B_q}=\sum_{r+q=n}{A_{r-1}\otimes B_q}=(\mathrm{S}A\otimes B)_n \ .$$
The verification that the differentials agree needs a little care: for $a\in A_{r-1}$ and $b\in B_q$, we have
\begin{eqnarray*}
d^{\mathrm{S}(A\otimes B)}(a\otimes b)=-d^{A\otimes B}(a\otimes b)=-(da\otimes b+(-1)^{r-1}a\otimes db) \ ,
\end{eqnarray*}
and
\begin{eqnarray*}
d^{\mathrm{S}A\otimes B}(a\otimes b)=-da\otimes b+(-1)^{r}a\otimes db \ .
\end{eqnarray*}
We leave the hom assertions to the reader.
\end{proof}

Let $f\dd A\to B$ be a chain map. The {\em mapping cone} $\mathrm{Mc}f$ of $f$ is the chain complex
defined by $(\mathrm{Mc}f)_n=B_n\oplus A_{n-1}$ and 
\begin{eqnarray}
d = \scriptsize{\left[ 
\begin{array}{ccc}
d & f \\
0 & -d 
\end{array} \right]} \ .
\end{eqnarray}
There is a short exact sequence
\begin{eqnarray}\label{sesMc}
0\lra B\xra{\scriptsize{\begin{bmatrix} 1  \\
0  
\end{bmatrix} }}\mathrm{Mc}f \xra{\scriptsize{\begin{bmatrix} 0 \ 1 
\end{bmatrix} }}\mathrm{S}A\lra 0 
\end{eqnarray}
of chain maps.
A little calculation with $2\times 2$-matrices proves:

\begin{proposition}\label{htpyMc}
If $f\dd A\to B$ is a chain map and 
$u\dd \mathrm{U}\mathrm{S}A \to \mathrm{U}B$ is a graded map
then 
\begin{eqnarray*}
\scriptsize{\begin{bmatrix} 1 & u \\
0 & 1 
\end{bmatrix}} 
\dd \mathrm{Mc}f \xra{\phantom{AAA}} \mathrm{Mc}(f+d(u))
\end{eqnarray*}
is an invertible chain map.
\end{proposition}

\begin{proposition}\label{absMc}
Given chain maps $i\dd B\to C$, $p\dd C\to \mathrm{S}A$ and graded maps
$j\dd \mathrm{U}\mathrm{S}A \to \mathrm{U}C$, $q\dd \mathrm{U}C\to \mathrm{U}B$
satisfying 
$$p\circ i=0 \ , \ q\circ i=1 \ , \ p\circ j = 1 \ , \ i\circ q+j\circ p=1 \ ,$$ 
it follows that $q\circ d(j) : A\to B$ is a chain map and 
\begin{eqnarray*}
[i,j]\dd \mathrm{Mc}(q\circ d(j))\lra C
\end{eqnarray*}
is an invertible chain map.
\end{proposition}
\begin{proof}
The equations imply that we have a short exact sequence
\begin{eqnarray*}
0\to B\xra{i}C\xra{p}\mathrm{S}A\to 0
\end{eqnarray*}
of chain maps since it is split exact in each dimension.
It follows that exactness is preserved by any additive functor.
Consequently,
\begin{eqnarray*}
0\to [\mathrm{S}A,B]\xra{[p,1]}[C,B]\xra{[i,1]}[B,B]\to 0
\end{eqnarray*}
is also a short exact sequence of chain maps. The corresponding
long exact homology sequence \cite{CarEil} includes the connecting 
graded morphism
\begin{eqnarray*}
\partial : \mathrm{H}[B,B]\lra \mathrm{H}[A,B]
\end{eqnarray*}
which takes the homology class of the cycle $1_B$ to the homology class
of the cycle $q\circ d(j)$. So $q\circ d(j) : A\to B$ is a chain map.
The remaining matrix calculations are left to the reader. 
\end{proof}

\begin{remark}\label{Mcascok}
{\em From the long exact homology sequence of \eqref{sesMc}, we see that the mapping cone of 
an identity chain map has zero homology. 
For any chain map $f\dd A\to B$, there is a short exact sequence
\begin{eqnarray*}
0\to A\xra{i'}B\oplus \mathrm{Mc}1_A\xra{p'}\mathrm{Mc}f\to 0 
\end{eqnarray*}
of chain maps where $i'= \scriptsize{{\left[\begin{array}{ccc}
-f \\
1_A \\
0 
\end{array}\right]}}$, $p'= \scriptsize{{\left[\begin{array}{ccc}
1_B & f & 0 \\
0 & 0 & 1
\end{array}\right]}}$.
If we take $j' = \scriptsize{{\left[\begin{array}{ccc}
1 &  0 \\
0 & 0 \\
0 & 1
\end{array}\right]}}$ and $q' = \scriptsize{{\left[\begin{array}{ccc}
0 &  1 & 0
\end{array}\right]}}$, we see that the equations of Proposition~\ref{absMc} hold.
The middle term in the short-exact sequence has the same homology as $B$ 
and $i'$ induces the same morphism on homology as $f$. So 
the construction can be used to replace the homology of an arbitrary chain map 
$f$ by the homology of an inclusion $i'$. (It is standard to use the mapping 
cylinder for that purpose.)}
\end{remark}

\begin{proposition}\label{LU=Mc1}
For any chain complex $A$, there is a natural isomorphism
$${\scriptsize{\left[ \begin{array}{ccc}
1 & d\\
-d & 1 
\end{array}\right]}}  \dd \mathrm{Mc}1_{\mathrm{S}^{-1}A} \cong \mathrm{LU}A \ .$$
Moreover, $\mathrm{L}\mathbb{Z}\otimes A \cong \mathrm{Mc}1_{\mathrm{S}^{-1}A}$
and $[\mathrm{L}\mathbb{Z},A] \cong \mathrm{Mc}1_A$.
There are corresponding results for $\mathrm{R}=\mathrm{LS}$. 
\end{proposition}
 
 \begin{proposition}
In the monoidal category $\mathrm{DGAb}$, there is a duality 
$$\mathrm{L}\mathbb{Z}\dashv\mathrm{R}\mathbb{Z} \ .$$
\end{proposition}

 \begin{proposition}\label{Urepble}
For each $n\in \mathbb{Z}$, the $\mathrm{Ab}$-functor $\mathrm{U}_n : \mathrm{DGAb}\lra \mathrm{Ab}$, 
obtained by following $\mathrm{U} : \mathrm{DGAb}\lra \mathrm{GAb}$ by evaluation at $n$, 
is represented by the chain complex $\mathrm{S}^n\mathrm{L}\mathbb{Z}$.
\end{proposition}
\begin{proof}
We have isomorphisms
\begin{eqnarray*}
\mathrm{DGAb}(\mathrm{S}^n\mathrm{L}\mathbb{Z},A) & \cong & \mathrm{DGAb}(\mathrm{L}\mathbb{Z},\mathrm{S}^{-n}A) \\ 
& \cong &  \mathrm{GAb}(\mathbb{Z},\mathrm{U}\mathrm{S}^{-n}A) \\
& \cong & (\mathrm{U}\mathrm{S}^{-n}A)_0 \ \cong \ \mathrm{U}_nA
 \end{eqnarray*}
which are natural in chain complexes $A$.
\end{proof}

Let $\CL$ denote the $\mathrm{Ab}$-category obtained as the full image of the functor $\mathbb{Z}\to \mathrm{DGAb}$
taking $n$ to $\mathrm{S}^n\mathrm{L}\mathbb{Z}$.
So the objects of $\CL$ are the integers and, using Proposition~\ref{Urepble}, the homs are
\begin{eqnarray*}
\CL(m,n) = \mathrm{DGAb}(\mathrm{S}^m\mathrm{L}\mathbb{Z},\mathrm{S}^n\mathrm{L}\mathbb{Z}) 
\cong (\mathrm{L}\mathbb{Z})_{m-n} = \left\{
\begin{array}{rl}
\mathbb{Z} & \text{for } n=m, m+1 \\
0 & \text{ otherwise } .
\end{array} \right.
\end{eqnarray*}  
The compositions $\CL(m,n)\ot \CL(\ell,m)\to \CL(\ell,n)$ are zero unless $m=n$ or $m=\ell$ in which
cases they are the canonical isomorphisms.  

 \begin{proposition}\label{denseness}
The singular $\mathrm{Ab}$-functor $\mathrm{DGAb} \to [\CL^{\mathrm{op}}, \mathrm{Ab}]$ for 
the inclusion $\CL\hookrightarrow \mathrm{DGAb}$ is an equivalence.
In particular, $\CL\hookrightarrow \mathrm{DGAb}$ is a dense $\mathrm{Ab}$-functor. 
\end{proposition}
\begin{proof}
The second sentence follows from the first since denseness of the inclusion precisely means that the singular
$\mathrm{Ab}$-functor is fully faithful \cite{KellyBook}. 

Consider the ordered set $\mathbb{Z}$ of integers as a category in the usual way and take the 
free pointed-set-enriched (``pointed'') category on it. Let $\mathbb{D}$ be the pointed category 
obtained by imposing the relation that the morphisms $m\to n$ for all $m\le n-2$ are
equal to the point $0$ in $\mathbb{D}(m,n)$. Straight from its definition, $\mathrm{DGAb}$ is equivalent
as a pointed category to the pointed functor category $[\mathbb{D}^{\mathrm{op}}, \mathrm{Ab}]_*$.
From the definition of $\CL$, we see that it is (isomorphic to) the free $\mathrm{Ab}$-category on the pointed category
$\mathbb{D}$ so that we have an equivalence $\mathrm{DGAb} \to [\CL^{\mathrm{op}}, \mathrm{Ab}]$ as pointed
categories and also that the equivalence underlies the singular $\mathrm{Ab}$-functor of the Proposition.
An $\mathrm{Ab}$-functor is an equivalence if and only if it is an equivalence as a (pointed) functor.  
\end{proof}

\begin{remark}\label{moreonL} 
{\em \begin{itemize}
\item[(i)] Under the equivalence of Proposition~\ref{denseness}, the chain complexes $\mathrm{S}^n\mathrm{L}\mathbb{Z}$
correspond to the representables in $[\CL^{\mathrm{op}}, \mathrm{Ab}]$. Compare Example 6.6 of \cite{PMcC2000}.
Indeed, for any additive category $\CA$, we have an equivalence $\mathrm{DG}\CA \simeq [\CL^{\mathrm{op}}, \CA]$ of additive categories. Also, the isomorphism $\mathrm{DG}\CA \cong \mathrm{DG}(\CA^{\mathrm{op}})^{\mathrm{op}}$ of \eqref{cochain}
is induced by an isomorphism $\CL^{\mathrm{op}} \cong \CL$.
\item[(ii)] The closed monoidal structure on $\mathrm{DGAb}$ transports across the equivalence to one on 
$[\CL^{\mathrm{op}}, \mathrm{Ab}]$. Using the theory of Day convolution \cite{DayConv}, we obtain a promonoidal structure on $\CL$
by looking at tensor products of representables.
The equation
\begin{eqnarray*}
\mathrm{L}\mathbb{Z}\otimes \mathrm{L}\mathbb{Z} \cong \mathrm{L}\mathbb{Z} \oplus \mathrm{S}^{-1}\mathrm{L}\mathbb{Z}
\end{eqnarray*}
implies 
\begin{eqnarray*}
\mathrm{S}^m\mathrm{L}\mathbb{Z}\otimes \mathrm{S}^n\mathrm{L}\mathbb{Z} \cong \mathrm{S}^{n+m}\mathrm{L}\mathbb{Z} \oplus \mathrm{S}^{n+m-1}\mathrm{L}\mathbb{Z} \ ,
\end{eqnarray*}   
showing that tensor products of representables are finite coproducts of representables. 
Therefore the induced promonoidal structure on the finite coproduct completion of the 
$\mathrm{Ab}$-category $\CL$ is nearly monoidal; all it lacks is the monoidal unit.
\end{itemize}}
\end{remark}

 \begin{proposition}\label{conicalgenerating}
For every chain complex $A$, there exists a reflexive coequalizer of chain maps of the form
\begin{eqnarray*}
\sum_{\theta\in \Theta}{\mathrm{S}^{m_\theta}\mathrm{L}\mathbb{Z}} \ \rightrightarrows \sum_{\phi\in \Phi}{\mathrm{S}^{n_\phi}\mathrm{L}\mathbb{Z}}\ra A \ . 
\end{eqnarray*}
 \end{proposition}
\begin{proof}
Recall that every $\mathrm{Ab}$-weighted colimit $\mathrm{colim}(J,K)$, for $\mathrm{Ab}$-functors 
$J : \CD^{\mathrm{op}}\to\mathrm{Ab}$ and  $K : \CD\to\CX$ with $\CX$ cocomplete, can be constructed as a 
reflexive coequalizer of the form
\begin{eqnarray*}
\sum_{\lambda\in \Lambda}{KD_{\lambda}} \ \rightrightarrows \sum_{\gamma\in \Gamma}{KE_{\gamma}}\ra \mathrm{colim}(J,K) \ . 
\end{eqnarray*}
in $\CX$. One way to see this is to take the usual (see \cite{KellyBook}) construction of the weighted colimit as a
coequalizer of reflexive pair of morphisms 
$$\sum_{D,E\in \CD}{(\CD(D,E)\otimes JE)\otimes KD} \ \rightrightarrows \sum_{D\in \CD}{JD\otimes KD} \ ,$$
then present the abelian groups $\CD(D,E)\otimes JE$ and $JD$ as reflexive coequalizers of morphisms
between free abelian groups. The tensor $\mathbb{Z}\Xi\otimes X$ of a free abelian group $\mathbb{Z}\Xi$
with $X\in \CX$ is a coproduct of $\Xi$ copies of $X$. 
So we have a $3\times 3$-diagram of reflexive coequalizers of morphisms between coproducts of objects in the image of $K$.
The main diagonal is therefore a reflexive coequalizer of a pair of morphisms; see the lemma on page 1 of \cite{JstneBE}. 

By Proposition~\ref{denseness}, our object $A\in \mathrm{DGAb}$ is an $\mathrm{Ab}$-weighted colimit of the dense inclusion 
$\CL\hookrightarrow \mathrm{DGAb}$ whose values are the complexes $\mathrm{S}^n\mathrm{L}\mathbb{Z}$.      
\end{proof}

A {\em double chain complex} in an additive category $\CA$ is a chain complex in $\mathrm{DG}\CA$.
We will use $\delta$ in place of $d$ for the differentials in $\mathrm{DG}\CA$.
So, a double complex $A$ consists of chain complexes $A_m$, for $m\in \mathbb{Z}$, with chain maps 
$\delta_m : A_m \to A_{m-1}$ satisfying $\delta_m\circ \delta_{m+1}=0$; we also have the
differentials $d^{A_m}_n :A_{m, n}\to A_{m, n-1}$. 

Suppose $A$ is a double chain complex in an additive category $\CA$ with countable coproducts. 
The {\em total complex} $\mathrm{Tot}A$ of $A$ is the chain complex with underlying object in  
$\mathrm{G}\CA$ defined by
\begin{eqnarray}\label{Tot}
\mathrm{Tot}A = \bigoplus_{m\in \mathbb{Z}} \mathrm{S}^mA_m
\end{eqnarray}
and with differential defined by the commutative diagram 
\begin{eqnarray*}
\xymatrix{
\mathrm{S}^mA_m \ar[rr]^-{\scriptsize{\begin{bmatrix} \mathrm{S}^m\delta_m  \\
d^{\mathrm{S}^mA_m}  
\end{bmatrix}}} \ar[d]_-{\mathrm{inj}_m}  &&  \mathrm{S}^mA_{m-1}\oplus \mathrm{S}^{m+1}A_m \ar[d]^-{\mathrm{inj}_{m-1,m}} \\
 \oplus_{m\in \mathbb{Z}} \mathrm{S}^mA_m \ar[rr]_-{d^{\mathrm{Tot}A}} &&  \oplus_{m\in \mathbb{Z}} \mathrm{S}^{m+1}A_m
  \ .}
\end{eqnarray*}
In the case $\CA = \mathrm{Ab}$, the formula for the differential becomes, for $a\in A_{m,n-m}$, 
\begin{eqnarray*}
d^{\mathrm{Tot}A}(a) = \delta_{m,n}(a) + (-1)^md^{A_m}_{n-m}(a) \in A_{m-1,n-m}\oplus A_{m,n-m-1} \ . 
\end{eqnarray*}

 \section{DG-Categories}\label{DG-cat}

A {\em DG-category} is a category $\CA$ with homs enriched in the monoidal category $\mathrm{DGAb}$ of chain complexes.
This was a motivating example for \cite{EilKel1966} and was developed somewhat in \cite{0,25}. 

To be explicit, $\CA$ has objects and, for objects $A$ and $B$, a chain
complex $\CA(A,B)$. The elements $u\dd A\to B$ of $\CA(A,B)_n$ are called {\em protomorphisms} of degree $n$.
The $n$-cycles of $\CA(A,B)$ are called {\em chain maps} of degree $n$. 
The ``underlying category'' $\CA_0$ of $\CA$ has
the same objects as $\CA$ and has chain maps of degree $0$ as morphisms. 
When we say ``chain map'' without mention the degree, we mean these morphisms of degree 0;
these are the morphisms of the underlying ordinary category of $\CA$.
Composition provides chain maps
\begin{eqnarray}
\circ \dd \CA(B,C)\ot \CA(A,B) \to \CA(A,C)
\end{eqnarray}
yielding, for protomorphisms $u\dd A\to B$ and $v\dd B\to C$ of degree $p$ and $q$
respectively, a protomorphism $v\circ u\dd A \to C$ of degree $p+q$ satisfying the equation
\begin{eqnarray}
d(v\circ u)=d(v)\circ u + (-1)^qv\circ d(u) \ .
\end{eqnarray}
In particular, a composite of chain maps is a chain map.
Composition is associative and each object $A$ has a chain map $1_A\dd A\to A$ of degree 0
which is an identity for composition. 
\begin{example}
{\em For any additive category $\CA$, the additive category $\mathrm{DG}\CA$ becomes a DG-category
by defining its hom chain complexes $\mathrm{DG}\CA [B,C]$ by
 \begin{eqnarray}
& \mathrm{DG}\CA [B,C]_n=\prod_{r-q=n}{\CA(B_q, C_r)} \nonumber \\
& (df)_q=d(f_q)-(-1)^nf_{q-1}d \ \text{ for } f = (f_q)_{q\in \mathbb{Z}}\in  \mathrm{DG}\CA [B,C]_n \ ,
\end{eqnarray}
abstracting \eqref{chainhom}.}
\end{example}

The {\em opposite} $\CA^{\mathrm{op}}$ of a DG-category $\CA$ is a DG-category with the same objects 
as $\CA$, with hom complexes $\CA^{\mathrm{op}}(A,B)=\CA(B,A)$, and where the composite of
protomorphisms $h\in \CA^{\mathrm{op}}(A,B)_m$, $k\in \CA^{\mathrm{op}}(B,C)_n$ in $\CA^{\mathrm{op}}$
is $(-1)^{mn}h\circ k \in \CA(C,A)_{m+n}$.  

The {\em tensor product} $\CA\otimes \CX$ of DG-categories $\CA$ and $\CX$ is the DG-category with
objects pairs $(A,X)$ of objects $A\in \CA$, $X\in \CX$, with hom complexes 
$$(\CA\otimes \CX)((A,X),(B,Y))= \CA(A,B)\otimes \CX(X,Y) \ ,$$ 
and composition defined by $(k,v)\circ (h,u)=(-1)^{qm}(k\circ h,v\circ u)$ where $m$, $q$ are
the degrees of $h$, $v$, respectively.

A {\em DG-functor} $F\dd \CA\to \CB$ between DG-categories assigns an object $FA$ of $\CB$
to each object $A$ of $\CA$ and a chain map $$F=F_{A,B}\dd \CA(A,A')\to \CB(FA,FA')$$ to each
pair of objects $A$, $B$, such that composition and identities are preserved.
These are precisely the functors enriched in the monoidal category $\mathrm{DGAb}$ 
in the sense of \cite{EilKel1966, KellyBook}. 

A {\em protonatural transformation (p.n.t)} $\theta\dd F\Ra G$ {\em of degree $n$} 
between DG-functors $F, G\dd \CA\to \CB$
is a family of protomorphisms $\theta_A\dd FA\to GA$ of degree $n$ such that
\begin{eqnarray*}
\theta_B\circ Ff=(-1)^{pn}Gf\circ \theta_A
\end{eqnarray*}
for all protomorphisms $f\dd A\to B$ of degree $p$ in $\CA$.
A {\em DG-natural transformation} is a p.n.t. $\theta$ for which each $\theta_A$ is a chain map of degree $0$;
these are precisely the natural transformations enriched in $\mathrm{DGAb}$ 
in the sense of \cite{EilKel1966, KellyBook}.

Let $\mathrm{DG\text{-}Cat}$ denote the 2-category of DG-categories, DG-functors, and DG-natural transformations.
It is a symmetric monoidal 2-category with respect to the tensor product $\CA\otimes \CB$. 

This gives the ingredients for an explicit description of the $\mathrm{DGAb}$-enriched functor category $[\CA,\CB]$
which we will call a {\em DG-functor DG-category}. The objects are the DG-functors $F\dd \CA\to \CB$.
The chain complex $[\CA,\CB](F,G)$ has p.n.t's $\theta$ of degree $n$ as elements in degree $n$ with differential
defined by $d(\theta)_A=d(\theta_A)$. (The reader can check that $d(\theta)$ is a p.n.t of degree $n-1$.)
As usual for categories enriched over a closed monoidal complete and cocomplete base, 
we have the isomorphism of DG-categories
\begin{eqnarray}
[\CA\ot \CB,\CC]\cong [\CA,[\CB,\CC]] \ .
\end{eqnarray}

Categories enriched in $\mathrm{GAb}$ are called {\em G-categories} (or ``graded categories'').  
By applying the strong-closed strong-monoidal functor $\mathrm{U}\dd \mathrm{DGAb}\to \mathrm{GAb}$ on hom complexes, we obtain a 
strong monoidal 2-functor $\mathrm{U}_*\dd \mathrm{DG\text{-}Cat}\to \mathrm{G\text{-}Cat}$.
Since the protonaturality condition does not involve the differentials and because $\mathrm{U}$ is strong closed, we have:

\begin{proposition}\label{protonat}
\begin{enumerate}
\item For DG-functors $F,G\dd \CA\to \CB$ between DG-categories, there is an identification
\begin{eqnarray*}
\mathrm{U}_*[\CA,\CB](F,G) = [\mathrm{U}_*\CA,\mathrm{U}_*\CB](\mathrm{U}_*F,\mathrm{U}_*G)
\end{eqnarray*}
in $\mathrm{GAb}$. 
\item For DG-functors $F,G\dd \CA\to \mathrm{DGAb}$ between DG-categories, there is an identification
\begin{eqnarray*}
\mathrm{U}_*[\CA,\mathrm{DGAb}](F,G) = [\mathrm{U}_*\CA,\mathrm{GAb}](\mathrm{U}F,\mathrm{U}G)
\end{eqnarray*}
in $\mathrm{GAb}$. 
\end{enumerate}
\end{proposition}

\section{Module-weighted colimits}

The theory of limits and colimits for categories enriched in a symmetric closed monoidal category is well developed
in the literature; see \cite{KellyBook, BorceuxBk2}.
We use the term ``weighted limit'' (as suggested in \cite{18}) rather than ``indexed limit''.

For a DG-category $\CC$, a {\em right DG-module} is a DG-functor 
$$M\dd \CC^{\mathrm{op}}\to \mathrm{DGAb} \ .$$
For any protomorphism $f\dd U\to V$ in $\CC$, we put $y\cdot f = (Mf)y$ for $y\in MV$.
For $f$ of degree $m$ and $g$ of degree $n$, notice then that
$$z\cdot (g\circ f)=(-1)^{nm}(z\cdot g)\cdot f \ .$$
  
The {\em $M$-weighted colimit} of a DG-functor $F\dd \CC\to \CA$ is an object 
$C = \mathrm{colim}(M,F)$ of $\CA$ equipped with chain isomorphisms
\begin{eqnarray*}
\pi_{F, A} : \CA(C,A)\cong [\CC^{\mathrm{op}},\mathrm{DGAb}](M,\CA(F-,A)) \ ,
\end{eqnarray*}
DG-natural in $A$. By Yoneda, the DG-natural transformation $\pi_F$ is determined by the DG-natural transformation
$\gamma_F = \pi_{F, C}(1_C) : M\to \CA(F-,C)$.
The $M$-weighted colimit $C = \mathrm{colim}(M,F)$ is {\em preserved} by the DG-functor $T : \CA\to \CX$ 
when the morphisms  
\begin{eqnarray*}
\CX(TC,X)\lra [\CC^{\mathrm{op}},\mathrm{DGAb}](M,\CA(TF-,X)) 
\end{eqnarray*}
of the DG-natural family indexed by $X$, corresponding to the
composite $$M\xra{\gamma_F} \CA(F-,C)\xra{T}\CX(TF-,TC) \ ,$$ are invertible. 
An $M$-weighted colimit in $\CA$ is {\em absolute} when it is preserved by all DG-functors out of $\CA$. 

An important special case is where $\CC = \CI$ is the DG-category with one object $0$
and with $\CI(0,0)=\mathbb{Z}$. A right $\CI$-module $M$ amounts to a chain complex $C$, 
while a DG-functor $F\dd \CI \to \CA$ amounts to an object $A$ of $\CA$.
Then $\mathrm{colim}(M,F) = C\ot A$, called the {\em tensor} of $C$ with $A$, satisfying
\begin{eqnarray*}
\CA(C\ot A,D)\cong \mathrm{DGAb}(C,\CA(A,D)) \ .
\end{eqnarray*}
  
This and Proposition~\ref{Sten} motivate the definition $\mathrm{S}A=\mathrm{S}\mathbb{Z}\ot A$
for {\em suspension} of an object $A$ in $\CA$. Clearly an object $B$ is isomorphic to $\mathrm{S}A$
when there exists an invertible chain map $A\to B$ of degree $-1$.

We call the DG-category $\CA$ {\em stable} when it admits tensoring with both $\mathrm{S}\mathbb{Z}$
and $\mathrm{S}^{-1}\mathbb{Z}$; then, for all integers $n$, we have 
$$\mathrm{S}^{n}A = \mathrm{S}^{n}\mathbb{Z}\ot A \ .$$
We call $\mathrm{S}^{-1}A$ the {\em desuspension} of $A$.  

Similarly, Proposition~\ref{LU=Mc1} motivates defining $\mathrm{Mc1_A}= \mathrm{L}\mathbb{Z}\ot \mathrm{S}A$
which we call the {\em mapping cone of the identity of $A$};
there is a canonical chain map $i_1\dd A\to \mathrm{Mc1_A}$ and a protomorphism
$q_1\dd \mathrm{Mc1_A}\to A$ with $q_1\circ i_1=1_A$.  

Conical limits are another important case of weighted colimits.
Examples are coproducts, coequalizers and cokernels.

Binary coproducts are of course {\em direct sums}: they also provide the binary product.
The direct sum $A\oplus B$ in $\CA$ is an object equipped with chain maps
$i\dd A\to A\oplus B$, $j\dd B\to A\oplus B$, $p\dd A\oplus B \to B$, $q\dd A\oplus B \to A$
satisfying the equations  
\begin{eqnarray}\label{directsum}
p\circ i=0 \ , \ q\circ i=1 \ , \ p\circ j = 1 \ , \ i\circ q+j\circ p=1 \ .
\end{eqnarray}
It follows that $q\circ j = 0$ showing the duality of the conditions.

Suspension, desuspension, finite direct sum, and mapping cone are all absolute colimits
since they can be expressed in terms of equations that are preserved by all DG-functors.
We now introduce another absolute colimit.

A {\em protosplitting} of a chain map $f:A\to B$ is a protomorphism $t:B\to A$ satisfying $f\circ t \circ f = f$.
Notice in this case $e = 1-f\circ t$ is an idempotent protomorphism and $e \circ f = 0$.

If $f$ is a chain monomorphism then the condition becomes $t \circ f = 1_A$ 
and we say $A$ is a {\em protosplit subobject} of $B$. 
If $f$ is a chain epimorphism then the condition becomes $f \circ t = 1_B$ 
and we say $B$ is a {\em protosplit quotient} of $A$. 
\begin{proposition}
Cokernels of protosplit chain maps are absolute.
\end{proposition}
\begin{proof}
Suppose $f:A\to B$ is protosplit by $t$.

Assume that the chain map $w : B \to C$ is a cokernel for $f$.
Since $e = 1-f\circ t$ satisfies $e\circ f=0$, there exists a unique protomorphism $s : C\to B$ with $s \circ w= e$.
Then $w\circ s\circ w = w\circ e = w - w\circ x\circ t = w$ implies $w\circ s = 1_C$.
So the chain map $w$ with $s$ form a protosplitting of the idempotent $e$. 
 
 Conversely, assume $e$ has a splitting $e = s \circ w$ with $w$ a chain map.
Then $w$ is a cokernel for $f$.   
To see this, suppose $g : B \to X$ is a protomorphism with $g\circ f=0$. 
Then $g\circ s \circ w = g\circ e = g - g\circ f\circ t = g$ and, if $h\circ w = g$ then $h = g\circ s$.
\end{proof}

Another class of absolute colimit is coequalizers of {\em protosplit parallel pairs} 
$u, v\dd A\to B$ of chain maps.
Such a pair requires that there should be a protomorphism $t\dd B\to A$ such that 
$$u\circ t = 1_B \ , \ v\circ t\circ u=v\circ t\circ v \ .$$

\begin{proposition}
Suppose the DG-category $\CA$ admits direct sums. 
In $\CA$, the existence of coequalizers of protosplit parallel pairs of chain maps is 
equivalent to the existence of cokernels of protosplit chain maps.
\end{proposition}
\begin{proof}
If $(u,v)$ is a protosplit parallel pair of chain maps then it is easy to see that $u-v$ is a protosplit chain map for which a cokernel
is a coequalizer of the pair.
Conversely, suppose $f \dd A\to B$ is a chain map and $t\dd B\to A$ a protomorphism with $f\circ t\circ f=f$.
The chain maps $u=\begin{bmatrix} 0 & 1 \end{bmatrix}$ and $v= \begin{bmatrix} f & 1 \end{bmatrix}$
from $A\oplus B$ to $B$
are protosplit by $t=\scriptsize{\begin{bmatrix} -t \\ 1 \end{bmatrix}}$.
A coequalizer for $u,v$ is a cokernel for $f$. 
\end{proof}

\begin{proposition}\label{idemsplitauto}
If a DG-category $\CA$ admits coequalizers of protosplit parallel pairs of chain maps then chain idempotents split.
\end{proposition}
\begin{proof}
For an idempotent $e\dd A\to A$, take $u=t=1_A$ and $v=e$. This gives a (proto)split
parallel pair whose coequalizer splits $e$.
\end{proof}

In Section 1 Example 4 of \cite{25}, it was shown how to obtain mapping cone 
as a weighted limit. Here we shall produce it dually from weighted colimits.
Motivated by Remark~\ref{Mcascok}, we define the {\em mapping cone}
of a chain map $f\dd A\to B$ in $\CA$ to be the cokernel of the chain map
\begin{eqnarray}
i= \scriptsize{{\left[\begin{array}{ccc}
-f \\
i_1 
\end{array}\right]}} \dd A \xra{\phantom{ooo}} B\oplus \mathrm{Mc}1_A \ ,
\end{eqnarray}
which is indeed a protosplit monomorphism.
What Remark~\ref{Mcascok} implies is:

\begin{proposition}
If a DG-category $\CA$ admits direct sums, suspension, cokernels of protosplit chain monomorphisms and mapping cones of identities then it admits mapping cones.
\end{proposition}

\begin{example}\label{weight_for_tot}
{\em An example of a non-absolute weighted colimit is the total complex introduced in Section~\ref{Cocc}.
The formula there is for the total complex of a double chain complex $A$ which we can think of as a DG-functor
from the additive category $\CL^{\mathrm{op}}$, regarded as a DG-category with homs all put in degree zero,
to $\mathrm{DGAb}$. Now take any DG-category $\CA$. For any DG-functor $A : \CL^{\mathrm{op}}\to \CA$,
we define $\mathrm{tot}A = \mathrm{colim}(J,A)$ where the weight $J : \CL \to \mathrm{DGAb}$
is described as follows. It amounts to the cochain complex $J$ in $\mathrm{DGAb}$ defined by 
$$J^m = \mathrm{S}^m\mathrm{L}\mathbb{Z}, \ \text{with codifferentials} \ 
\delta^m = (\mathrm{S}^m\mathrm{L}\mathbb{Z}\xra{\mathrm{S}^m d^{\mathrm{L}\mathbb{Z}}} \mathrm{S}^{m+1}\mathrm{L}\mathbb{Z})  \ .$$
It is an exercise for the reader to check that, in the case $\CA = \mathrm{DGAb}$, 
we have $\mathrm{colim}(J,A) \cong \mathrm{Tot}A$ as in \eqref{Tot}. 
We offer the following diagram and point to Proposition~\ref{Urepble} as hints.
\begin{equation*}
\xymatrix{
J^m_p \ar[rr]^-{\theta^m} \ar[d]_-{\delta^m} && [A_m,B]_{p+n} \ar[d]^-{[\delta_{m+1},1]} \\
J^{m+1}_p \ar[rr]_-{\theta^{m+1}} && [A_{m+1},B]_{p+n}}
\end{equation*}
}
\end{example}

For a DG-category $\CC$, a {\em left DG-module} is a DG-functor $N\dd \CC\to \mathrm{DGAb}$. 
Given a right DG-module $M$ and a left DG-module $N$ for $\CC$, 
we define their {\em tensor product over $\CC$} by 
\begin{eqnarray}
M\ot_{\CC}N= \mathrm{colim}(M,N) = \int^U{MU\ot NU} = \sum_U{MU\ot NU}/\sim \ ,
\end{eqnarray}
where the congruence $\sim$ is generated by 
\begin{eqnarray*}
(y\cdot f)\ot x \sim y\ot (f\cdot x) \ \text{ in } \ (\sum_W{MW\ot NW})_n
\end{eqnarray*}
for $x\in (NU)_r$, $f\in \CC(U,V)_s$, $y\in (MV)_t$ with $r+s+t=n$. 

In terminology not exactly used in \cite{LawMetric} but inspired by it, and used in \cite{21} and elsewhere, 
we define a right DG-module $M$
over $\CC$ to be {\em Cauchy} when there exists a left DG-module $N$, a chain map
$$\eta\dd \mathbb{Z}\lra M\ot_{\CC}N$$
and a DG-natural transformation
$$\varepsilon_{U,V}\dd NU\ot MV\lra \CC(V,U)$$
such that diagram \eqref{snake} commutes.
\begin{eqnarray}\label{snake}
\begin{aligned}
\xymatrix{
MV \ar[rr]^-{1_{MV}} \ar[d]_-{\cong} &  & MV \ar[d]^-{\cong} \\
\mathbb{Z}\ot MV \ar[rd]_-{\eta\ot 1_{MV}} &  & {\int^U MU\ot \CC (V,U)} \\
& {\int^U MU\ot NU \ot MV} \ar[ru]_-{\phantom{aaa}{\int^U1_{MU}\ot \varepsilon_{U,V}}} &}
\end{aligned}
\end{eqnarray}
Equivalently (see Section 5.5 of \cite{KellyBook} for example), $M$ is Cauchy when it is ``small projective''
in $[\CC^{\mathrm{op}},\mathrm{DGAb}]$; that is, when the DG-functor 
$$[\CC^{\mathrm{op}},\mathrm{DGAb}](M,-)\dd [\CC^{\mathrm{op}},\mathrm{DGAb}] \lra \mathrm{DGAb}$$
represented by $M$ preserves weighted colimits. 

A right DG-module $M$ over $\CC$ is called {\em representable} or {\em convergent} when there
is an object $K$ of $\CC$ and a DG-natural isomorphism $\CC(-,K)\cong M$.
Convergent modules are Cauchy.
In agreement with terminology of \cite{LawMetric}, we call a DG-category $\CC$ {\em Cauchy complete} when every Cauchy right DG-module
over it is convergent.
The {\em Cauchy completion} $\mathrm{Q}\CC$ of a DG-category $\CC$ is the full sub-DG-category
of $[\CC^{\mathrm{op}},\mathrm{DGAb}]$ consisting of the Cauchy right $\CC$-modules; compare \cite{21}.

The Morita Theorem in the enriched context here gives:

\begin{theorem}
For any DG-categories $\CC$ and $\CD$, there is an equivalence 
$$[\CC^{\mathrm{op}},\mathrm{DGAb}]\simeq [\CD^{\mathrm{op}},\mathrm{DGAb}]$$
if and only if there is an equivalence $$\mathrm{Q}\CC\simeq \mathrm{Q}\CD \ ,$$
and this, if and only if there is an equivalence
$$[\CC,\mathrm{DGAb}]\simeq [\CD,\mathrm{DGAb}] \ .$$
Moreover, $\mathrm{QQ}\CC\simeq \mathrm{Q}\CC$, so
$$[\mathrm{Q}\CC,\mathrm{DGAb}]\simeq [\CC,\mathrm{DGAb}] \ .$$    
\end{theorem}

A weighted colimit in a DG-category $\CA$ is called {\em absolute} when it is preserved by all DG-functors out of $\CA$.

\begin{theorem}[\cite{23}]
A right $\CC$-module $M$ is Cauchy if and only if all $M$-weighted colimits are absolute.
\end{theorem} 

From \cite{23}, we obtain:

\begin{proposition}
Every Cauchy complete DG-category admits cokernels of protosplit chain maps.
\end{proposition}  

\begin{corollary}\label{oneway}
The right DG-modules providing the weights for direct sums, suspension, desuspension, 
mapping cones of identities, and cokernels of protosplit chain maps are all Cauchy.
\end{corollary}
\begin{proof}
All these colimits are determined by equations preserved by DG-functors; so they are absolute.
\end{proof}

\section{Cauchy G-modules}\label{CcGm}

The inclusion $\mathrm{GAb}\hookrightarrow \mathrm{DGAb}$ of graded abelian groups in chain complexes
(by taking $d=0$) is symmetric strong monoidal and strong closed.
So every G-category $\CX$ (that is, category enriched in $\mathrm{GAb}$) can be regarded as a DG-category. 
The distinction between chain maps and protomorphism disappears since every protomorphism
of a given degree is automatically a chain map of that degree; we simply call them morphisms of that degree.
The underlying category $\CX_0$ of $\CX$ has morphisms those of degree $0$. 

Cocompleteness for a G-category (in the $\mathrm{GAb}$-enriched sense) 
differs from cocompleteness as a DG-category (in the $\mathrm{DGAb}$-enriched sense).
Even the base $\mathrm{GAb}$ is not DG-complete.
It does not admit mapping cones of identities, for example. 
That is to be expected since the weight $\mathrm{L}\mathbb{Z}$ is not in $\mathrm{GAb}$. 

A {\em right G-module} over a G-category $\CX$ is a G-functor $M\dd \CX^{\mathrm{op}}\to \mathrm{GAb}$.  
Tensoring with graded abelian groups, suspension, desuspension, direct sum, coequalizers and cokernels are all
examples of colimits weighted by right G-modules. In particular, we may speak of {\em stable} G-categories.

A right G-module is {\em Cauchy} when the $N$ in the definition of Cauchy DG-module lands in $\mathrm{GAb}$.
However, this is the same as being Cauchy as a DG-module. 
For, forgetting the differentials in the values of any $N$ will provide an 
$N$ satisfying \eqref{snake} for the G-module. 

\begin{proposition}\label{GCauchymod}
A right G-module $M$ over a G-category $\CX$ is Cauchy if and only if there exist objects $E_1, \dots , E_n\in \CX$
and integers $m_1, \dots , m_n$ such that $M$ is a retract of 
$\bigoplus_{i=1}^{n}{\mathrm{S}^{m_i}\CX(-,E_i)}$ in $[\CX^{\mathrm{op}},\mathrm{GAb}]$.
\end{proposition}
\begin{proof}
To prove ``if'' simply observe that all the colimits used to obtain $M$ from the representables $\CX(-,E_i)$ are absolute. 

Conversely, let $M$ be a Cauchy right G-module over $\CX$.
There is a left G-module $N$ over $\CX$, a graded morphism
$\eta\dd \mathbb{Z}\lra M\ot_{\CX}N$
and a G-natural transformation
$\varepsilon_{U,V}\dd NU\ot MV\lra \CX(V,U)$
such that diagram \eqref{snake} commutes. 
Now $\eta(1)$ is represented by an element of degree $0$ in the coproduct $\sum_{X\in \CX}MX\ot NX$
of graded abelian groups; that is, there are objects $E_1,\dots , E_n$ of $\CX$ and elements
$x_i\in (ME_i)_{m_i}$, $y_i\in (NE_i)_{-m_i}$, $i=1,\dots , n$, with 
$$\eta(1)=\left[ \sum_{i=1}^nx_i\ot y_i \right] \ .$$ 
Condition \eqref{snake} becomes
\begin{eqnarray}\label{evsnake}
\sum_{i=1}^nx_i\cdot \varepsilon_{E_i,X}(y_i\ot u) = u
\end{eqnarray}
for all $X\in \CX$, $u\in (MX)_r$, $r\in \mathbb{Z}$.
Define $\tau_i\dd M\Lra \mathrm{S}^{m_i}\CX(-,E_i)$ to be the G-natural transformation
with components
\begin{eqnarray*}
(MX)_r\xra{y_i\ot 1}(NE_i)_{-m_i}\ot (MX)_r\xra{\varepsilon_{E_i,X}}\CX(X,E_i)_{r-m_i} \ .
\end{eqnarray*}
Let $\widehat{x_i}\dd \mathrm{S}^{m_i}\CX(-,E_i)\Lra M$ denote the G-natural transformation
corresponding under Yoneda to $x_i\in (ME)_{m_i}$; explicitly, 
$$\widehat{x_i}(X\xra{f}E_i)=x_i\cdot f \ .$$
Condition \eqref{evsnake} amounts to commutativity of the triangle
\begin{equation*}
\xymatrix{
M \ar[rd]_{1_M}\ar[rr]^{\tau \phantom{AAAAAA}}   && {\bigoplus_{i=1}^n\mathrm{S}^{m_i}\CX(-,E_i)} \ar[ld]^{\widehat{x}} \\
& M  &
}
\end{equation*} 
where $\tau = {\scriptsize{\left[\begin{array}{ccc}
\tau_1 \\
\vdots \\
\tau_n
\end{array} \right]}}$ 
and $\widehat{x} = {\scriptsize{\left[\begin{array}{ccc}
\widehat{x_1} & \dots & \widehat{x_n}
\end{array} \right]}}$.
So $M$ is a retract as required.
\end{proof}

\begin{proposition}\label{Gcase}
If a G-category $\CX$ is stable and admits finite direct sums then every Cauchy right $G$-module over $\CX$ is a retract of a convergent right $G$-module.
A G-category $\CX$ is Cauchy complete if and only if it is stable, and admits idempotent splittings and finite direct sums.
\end{proposition}
\begin{proof}
To see ``only if'' in the second sentence, we note that splitting of idempotents, finite direct sums, suspension and desuspension are absolute colimits, so any Cauchy complete G-category must admit them.

Conversely, and for the first sentence, suppose the stable G-category $\CX$ admits finite direct sums.
Let $M$ be a Cauchy right G-module over $\CX$; apply Proposition~\ref{GCauchymod} and use that notation. 
Now there is a G-natural isomorphism 
$$\bigoplus_{i=1}^n\mathrm{S}^{m_i}\CX(-,E_i) \cong \CX(-,\bigoplus_{i=1}^n\mathrm{S}^{m_i}E_i) \ .$$
Using Yoneda, we see that the idempotent of Proposition~\ref{GCauchymod} on the left hand side transports
to an idempotent $\CX(-,e)$ on the right hand side of the last isomorphism.
This proves the first sentence of the Proposition.
When it exists, let $K$ be a splitting of the idempotent $e$ on $\bigoplus_{i=1}^n\mathrm{S}^{m_i}E_i$. 
By Yoneda again, $M\cong \CX(-,K)$. So $M$ converges, proving ``if'' in the second sentence.  
\end{proof}

\begin{corollary}
The Cauchy completion $\mathrm{Q}\CX$ of a G-category $\CX$ is the closure 
of the representables in $[\CX^{\mathrm{op}},\mathrm{GAb}]$ under 
suspension and desuspension, direct sums, and splitting of idempotents.
\end{corollary}
 

\section{Cauchy DG-modules}\label{CcDGm}

The goal of this section is to prove our main result which is the converse of Corollary~\ref{oneway}.

\begin{proposition}\label{DGconical}
Suppose $\CC$ is a stable DG-category which admits tensoring by $\mathrm{L}\mathbb{Z}$. 
Every right DG-module $M$ over $\CC$ is a cokernel 
\begin{eqnarray}\label{cokerforM}
\sum_{i\in I}\CC(-,A_i) \xra{\phi} \sum_{j\in J}\CC(-,B_j) \xra{\gamma} M \ra 0
\end{eqnarray}
of a chain map between coproducts of representables in $[\CC^{\mathrm{op}},\mathrm{DGAb}]$.
\end{proposition}
\begin{proof}
From \cite{KellyBook} for example, we know that $M \in [\CC^{\mathrm{op}},\mathrm{DGAb}]$ is a coequalizer of
a reflexive pair of chain maps
\begin{eqnarray*}
\sum_{D,E\in \CC}{ME\otimes \CC(D,E)\otimes \CC(-,D)} \ \rightrightarrows \sum_{D\in \CC}{MD\otimes \CC(-,D)} \ra M
\end{eqnarray*}
We now use Proposition~\ref{conicalgenerating} to present the chain complexes $ME\otimes \CC(D,E)$ and $MD$
as reflexive coequalizers of coproducts of objects of the form $\mathrm{S}^n\mathrm{L}\mathbb{Z}$.
So we have a $3\times 3$-diagram of reflexive coequalizers of morphisms between coproducts of objects of the form
$\mathrm{S}^n\mathrm{L}\mathbb{Z}\otimes \CC(-,D)$ and such objects are isomorphic to representables $\CC(-,\mathrm{S}^n\mathrm{L}\mathbb{Z}\otimes D)$ because $\CC$ is stable and admits tensoring by $\mathrm{L}\mathbb{Z}$. 
Making these replacements in the $3\times 3$-diagram, we take the main diagonal to obtain a reflexive pair whose
coequalizer is $M$ \cite{JstneBE}. Then the required $\phi$ is the difference of the chain maps in that reflexive pair. \end{proof}

\begin{theorem}\label{DGcase}
If a stable DG-category $\CC$ admits finite direct sums and tensoring by $\mathrm{L}\mathbb{Z}$ 
then every Cauchy right DG-module over $\CC$ is a protosplit quotient of a representable DG-module.
The DG-category $\CC$ is Cauchy complete if and only if, furthermore, it admits cokernels of protosplit chain maps.
\end{theorem}
\begin{proof} 
Let $\CC$ be a stable DG-category with finite direct sums and tensoring by $\mathrm{L}\mathbb{Z}$. 

Assume $M\in [\CC^{\mathrm{op}},\mathrm{DGAb}] = : \CP \CC$ is Cauchy.
By Proposition~\ref{DGconical}, there is a cokernel diagram \eqref{cokerforM}
in the additive category of chain maps in $\CP \CC$.  

For each $i\in I$, there exists a finite subset $J_i$ of $J$ such that
\begin{eqnarray}
\phi_{A_i}(1_{A_i}) = \Sigma_{j\in J_i}x_{i j} \ .
\end{eqnarray}
Since $\phi_{A_i}$ is a chain map and $1_{A_i}$ is a 0-cycle in $\CC(A_i,A_i)$,
it follows that $\sum_{j\in J_i}x_{i j}\in \bigoplus_{j\in J_i}\CC(A_i,B_j)$ is a 0-cycle.
So each $x_{i j} :A_i \to B_j$ is a chain map in $\CC$.
For convenience, put $x_{i j} = 0$ for $j\notin J_i$.
By Yoneda, 
\begin{eqnarray}
\phi_X(X\xra{u}A_i) = \Sigma_{j\in J}x_{i j}u \ .
\end{eqnarray}

We also have a 0-cycle $m_j = \gamma_{B_j}(1_{B_j})\in MB_j$ for each $j\in J$.
By Yoneda, 
\begin{eqnarray}
\gamma_X(X\xra{v} B_j) = (Mv)m_j \ .
\end{eqnarray}

Since $M$ is Cauchy, $\CP \CC(M,-)$ preserves coproducts and cokernels, so
we have a regular epimorphism
\begin{eqnarray}\label{regepi}
\sum_{j\in J}\CP \CC(M,\CC(-,B_j)) \twoheadrightarrow \CP \CC(M,M)
\end{eqnarray}
in $\mathrm{DGAb}_0$ induced by $\CP \CC(M,\gamma)$. 
In particular, the identity endomorphism of $M$ is in the image of \eqref{regepi}.
So there exists a finite subset $J'$ of $J$ and protonatural transformations
\begin{eqnarray}
\sigma_j : M \lra \CC(-,B_j) \ , \ j\in J' \ ,
\end{eqnarray}
such that each $m\in MX$ can be written as
\begin{eqnarray}\label{em}
m = \Sigma_{j\in J'}M(\sigma_{j X}(m)) m_j \ .
\end{eqnarray}

Put $B = \sum_{j\in J'}B_j \in \CC$ and define $\sigma$ and $\gamma'$ by
the commutative diagram
\begin{eqnarray}
\begin{aligned}
\xymatrix{
M \ar[rr]^-{\sigma} \ar[rrd]_-{[\sigma_j]} && \CC(-,B) \ar[rr]^-{\gamma'} \ar[d]_-{\cong} && M  \\
 & & \Sigma_{j\in J'}{\CC (-,B_j)} \ar[rr]_-{\mathrm{incl.}} && \Sigma_{j\in J}{\CC (-,B_j)} \ . \ar[u]_-{\gamma}} 
\end{aligned}
\end{eqnarray}
Equation \eqref{em} is equivalent to
\begin{eqnarray}\label{retract}
\gamma' \circ \sigma = 1_M \ ,
\end{eqnarray}
which yields the first sentence of the theorem.

Put $e = \sigma_B\gamma '_B(1_B) \in \CC(B,B)$. By Yoneda,
\begin{eqnarray}\label{idemp_e_def}
\sigma\circ \gamma' = \CC(-,e) : \CC(-,B)\lra \CC(-,B) \ .
\end{eqnarray}
By \eqref{retract}, $e$ is idempotent and
\begin{eqnarray}
\gamma '_B(e) = \gamma '_B \sigma_B \gamma '_B(1_B) = \gamma '_B(1_B) \ .
\end{eqnarray}

It follows that 
\begin{eqnarray*}
\Sigma_{j\in J'}\mathrm{pr}_j (1_B- e) \in \mathrm{ker}\gamma_B = \mathrm{im}\phi_B \ .
\end{eqnarray*}
So there exists a finite subset $I'$ of $I$ and protomorphisms $t_i : B \to A_i$ for $i\in I'$
such that 
\begin{eqnarray}\label{1-e}
\phi_B(\Sigma_{i\in I'}t_i) = \Sigma_{i\in I'}\Sigma_{j\in J}x_{i j}t_i = \Sigma_{j\in J'}\mathrm{pr}_j(1-e) \ .
\end{eqnarray}
Put $A = \sum_{i\in I'}A_i$ in $\CC$ so that the $t_i$ define a protomorphism $t : B\to A$.
Then the $x_{i j}$ determine a chain map $x : A \to B$. The equation $\gamma \phi = 0$
restricts to $\gamma'\circ \CC(-,x) = 0$; from \eqref{idemp_e_def} and Yoneda we deduce $ex=0$. 
Also \eqref{1-e} becomes $x t = 1-e$. This tells us that $x$ is protosplit by $t$.     
\end{proof}

\begin{remark} Cokernels of protosplit chain maps do not follow from stability, finite direct sums, 
tensoring by $\mathrm{L}\mathbb{Z}$, and splitting of (chain) idempotents. {\em To see this, let $\CC$ be the
smallest full sub-DG-category of $\mathrm{DGAb}$ containing $\mathrm{L}\mathbb{Z}$, and closed under finite direct sums, suspensions, desuspensions and retracts. Using the Remark~\ref{moreonL} (ii) and the fact that the only retracts of
$\mathrm{L}\mathbb{Z}$ are $0$ and itself, we see that the objects of $\CC$ are those of the form 
$$\bigoplus_{i=1}^n\mathrm{S}^{n_i}\mathrm{L}\mathbb{Z} \ .$$
Yet the chain map $\mathrm{S}^{-1}\mathrm{L}\mathbb{Z}\xra{f} \mathrm{L}\mathbb{Z}$ with $f_{-1}= 1_{\mathbb{Z}}$ is protosplit by $\mathrm{L}\mathbb{Z}\xra{t}  \mathrm{S}^{-1}\mathrm{L}\mathbb{Z}$ with $t_{-1}= 1_{\mathbb{Z}}$, whereas the cokernel of $f$ is $\mathbb{Z}\notin \CC$.   
}
\end{remark}

Here is the DG-analogue of Proposition~\ref{GCauchymod}.

\begin{proposition}\label{DGCauchymod}
A right DG-module $M$ over a DG-category $\CC$ is Cauchy if and only if there exist objects $C_{i}, C'_{j} \in \CC$
and integers $m_i$, $n_j$ for $1\le i\le r$, $1\le j\le s$ such that $M$ is a protosplit quotient of the DG-module 
\begin{eqnarray*}
 \bigoplus_{i=1}^{r}{\mathrm{S}^{m_i}\CC(-,C_{i})} \oplus  \bigoplus_{j=1}^{s}{\mathrm{S}^{n_j}\mathrm{LU}\CC(-,C'_{j})} \ .
\end{eqnarray*}
\end{proposition}
\begin{proof}
A protosplit quotient $M$ of the displayed right $\CC$-module is an absolute colimit of Cauchy modules and so Cauchy.

Conversely, assume $M$ is Cauchy and let $J : \CC\to \CC'$ be the inclusion of $\CC$ into its free cocompletion $\CC'$ under
finite sums, suspensions, desuspensions, and tensoring by $\mathrm{L}\mathbb{Z}$. 
The extension 
$M' : \CC'^{\mathrm{op}} \to \mathrm{DGAb}$ of $M : \CC^{\mathrm{op}} \to \mathrm{DGAb}$ is the left Kan extension of
$M$ along $J$; see \cite{KellyBook}. It follows that 
$[\CC'^{\mathrm{op}},\mathrm{DGAb}](M',-) \cong [\CC^{\mathrm{op}},\mathrm{DGAb}](M,-J)]$
preserves colimits. So $M'$ is Cauchy with $\CC'$ satisfying the assumptions of the first sentence of Theorem~\ref{DGcase}; so $M'$
is a protosplit quotient of a representable $\CC'(-,D)$.
Therefore $M \cong M'\circ J^{\mathrm{op}}$ is a protosplit quotient of $\CC'(J-,D)$. 
Using Remark~\ref{moreonL}, objects of $\CC'$ are of the form
\begin{eqnarray*}
D\cong  \bigoplus_{i=1}^{r}{\mathrm{S}^{m_i}C_{i}} \oplus  \bigoplus_{j=1}^{s}{\mathrm{L}\mathbb{Z}\otimes \mathrm{S}^{n_j}C'_{j}} \ .
\end{eqnarray*}
Since $J$ is fully faithful and $\mathrm{L}\mathbb{Z}\otimes -\cong \mathrm{LU}$, the result follows.  
\end{proof}

\section{Complexes in DG-categories}

Let $\CA$ be a DG-category. The underlying additive category $\mathrm{Z}_{0 *}\CA$ of $\CA$ 
is obtained by keeping the same objects as $\CA$ and applying the $0$-cycle functor $\mathrm{Z}_{0}$ (see \eqref{Zed}) on
the hom chain complexes.
A {\em complex} in $\CA$ is a chain complex in $\mathrm{Z}_{0 *}\CA$ as per \eqref{chaincomplex};
that is, it consists of a family $A=(A_m)_{m\in\mathbb{Z}}$
of objects of $\CA$ together with chain maps
\begin{eqnarray}
\delta_{m}\dd A_m\lra A_{m-1} \ \text{ for all } \ m\in \mathbb{Z}
\end{eqnarray}
of degree of 0 subject to the equation $\delta_m\circ \delta_{m+1}=0$.

\begin{example}
{\em 
\begin{itemize}
\item[(i)] If $\CA$ is a mere additive category then a complex in $\CA$ is a chain complex in $\CA$.
\item[(ii)] For an additive category $\CA$, a complex in $\mathrm{DG}\CA$ is a double chain complex in $\CA$ as defined in Section~\ref{Cocc}.
\item[(iii)]  If $B\in \mathrm{DGAb}$, we can define a complex $\bar{B}$ in $\mathrm{DGAb}$ to be
the family of chain complexes $\mathrm{S}^{-m}B$, $m \in \mathbb{Z}$, equipped with the graded morphisms
$\delta_{m}\dd \mathrm{S}^{-m}B\to \mathrm{S}^{-m+1}B$ having $n$-th component $d^B : B_{n+m}\to B_{n+m-1}$.  
\end{itemize}}
\end{example}

For each DG-category $\CA$, we define a DG-category $\mathrm{DG}\CA$ whose objects
are the complexes in $\CA$. The hom chain complex $\mathrm{DG}\CA(A,B)$ as a graded
abelian group is defined to be equal to 
\begin{eqnarray*}
\prod_q\bigoplus_p\mathrm{S}^{p-q}\CA(A_q,B_p) \ ,
\end{eqnarray*}
so its elements of degree $n$ consist of those families $f = (f_{p, q}\dd A_q\to B_p)$ 
of protomorphisms of degree $n-p+q$ such that, 
for each $q \in \mathbb{Z}$, 
$f_{p q} = 0$ for all but a finite number of $p$;
the differential is defined by
\begin{eqnarray}
d(f)_{p,q} = (-1)^p d(f_{p,q}) + \left( \delta_p \circ f_{p+1,q} - (-1)^n f_{p,q-1} \circ \delta_{q} \right) \ .
\end{eqnarray}
Composition 
$$\circ \dd \mathrm{DG}\CA(B,C)_u\ot \mathrm{DG}\CA(A,B)_v \lra \mathrm{DG}\CA(A,B)_{u+v}$$
is defined by the formula
\begin{eqnarray}\label{complcomp}
(g\circ f)_{p q}=\sum_{r\in \mathbb{Z}} g_{p r}\circ f_{r q} \ .
\end{eqnarray}
The identity $1_A\in \mathrm{DG}\CA(A,A)_0$ is given by the Kronecker delta 
\begin{eqnarray*}
(1_A)_{p q} =
\begin{cases}
1_{A_p} & \text{if } p=q \\
0 & \text{ otherwise} \ .
\end{cases}
\end{eqnarray*}

\begin{example}
{\em If $\CA$ is a mere additive category then $\mathrm{DG}\CA$ is the DG-category $\mathrm{DG}\CA$ of chain complexes in $\CA$.}
\end{example}

\begin{proposition}
For each DG-category $\CA$ with a zero object $0$, the inclusion $\mathrm{i} : \CA \to \mathrm{DG}\CA$, where $\mathrm{i}A$ is $0$ for all degrees except degree 0 for which it is $A$, is dense.  
That is, the singular DG-functor
$$\widetilde{\mathrm{i}} : \mathrm{DG}\CA \lra [\CA^{\mathrm{op}},\mathrm{DGAb}] \ ,$$
taking $B \in \mathrm{DG}\CA$ to $\mathrm{Tot}\CA(-,B) \dd \CA^{\mathrm{op}}\to \mathrm{DGAb}$, is fully faithful.
\end{proposition}
\begin{proof}
We only provide an outline. Before anything else, we should see that the singular DG-functor is as claimed.
For $A\in \CA$, the definition of $\mathrm{DG}\CA$ gives 
$\mathrm{DG}\CA(\mathrm{i}A,B)_n = \bigoplus_q\CA(A,B_q)_{n-q}$.
So $$\mathrm{DG}\CA(\mathrm{i}A,B) = \bigoplus_q\mathrm{S}^q\CA(A,B_q) = \mathrm{Tot}\CA(A,B) \ .$$ 

Now we address the full faithfulness.
Take a protonatural family $\theta_A : \mathrm{Tot}\CA(A,B) \to \mathrm{Tot}\CA(A,C)$ of degree $n$ which composes
with the $q$-injection to give graded morphisms 
$$\theta_A^q : \mathrm{S}^q\CA(A,B_q)\to \bigoplus_p\mathrm{S}^p\CA(A,C_p) $$
of degree $n$.
Using the Yoneda Lemma for DG-categories, we have that these amount to elements 
$t_{q} \in \bigoplus_p\CA(B_q,C_p)_{q-p+n}$,
determining a unique $t \in \mathrm{DG}\CA(B,C)_n$ with $\widetilde{\mathrm{i}}(t) = \theta$. 
\end{proof}

\begin{proposition}\label{tot_as_ladj}
For any DG-category $\CA$ with zero object, the value at a complex $A$ in $\CA$ of a left adjoint to the 
DG-functor $\mathrm{i} : \CA \to \mathrm{DG}\CA$ is given by 
$\mathrm{tot}A = \mathrm{colim}(J,A)$ as defined in Example~\ref{weight_for_tot}.
In case $\CA = \mathrm{DG}\CC$ for an additive category $\CC$ with zero object, 
the DG-functor $\mathrm{i}$ has a left adjoint if and only if $\CC$ has countable coproducts.   
\end{proposition}
\begin{proof}
We have the graded isomorphisms
\begin{eqnarray*}
\mathrm{DG}\CA(A, \mathrm{i}X) & \cong & \prod_m\bigoplus_n\mathrm{S}^{n-m}\CA(A_m,\mathrm{i}X_n) \\
& \cong & \prod_m\mathrm{S}^{-m}\CA(A_m,X) \\
& \cong &\mathrm{lim}(J,\CA(A,X)) \\
& \cong & \CA (\mathrm{colim}(J,A),X)  
 \end{eqnarray*}
under which the differentials correspond.
This proves the first sentence. The second sentence uses the construction \eqref{Tot} of $\mathrm{tot}A$
to prove ``if''. For ``only if'', take a countable family $(X_n)_{n\in \mathbb{N}}$ of objects of $\CC$ and
look at the complex $A$ in $\mathrm{DG}\CC$ defined by $A_{m,n} = X_n$ when $0\leq n=-m$
and $A_{m,n} = 0$ otherwise, with all differentials zero morphisms; then 
$(\mathrm{tot}A)_0 \cong \bigoplus_{n\in \mathbb{N}}X_n$.   
\end{proof}

There are some obvious questions worthy of later consideration.
\begin{question}\begin{itemize}
\item[(a)] If $\CA$ is a Cauchy complete DG-category, is $\mathrm{DG}\CA$ also Cauchy complete?
\item[(b)] Is there a simple description of the completion of a DG-category $\CA$ with respect to $J$-weighted colimits? 
\end{itemize}
\end{question}

Here is an easy case of (a):

\begin{proposition}
If $\CA$ is a Cauchy complete additive category then $\mathrm{DG}\CA$ is a Cauchy complete DG-category. 
\end{proposition}
\begin{proof}
The construction of suspension and desuspension require nothing. Also finite direct sums and mapping cones of chain maps in $\mathrm{DG}\CA$ are straightforward as in $\mathrm{DGAb}$; they just use finite direct sums in $\CA$. 
Suppose we have a chain map $x : A \to B$ protosplit by $t : B\to A$. 
We have the idempotent protomorphism $e=1-xt : B\to B$ with 
$ex = 0$.
Since idempotents split in $\CA$ and thus in $\mathrm{G}\CA$, the chain map $x$ has a cokernel $w\dd B\to C$ in 
$\mathrm{G}\CA$. Then $wd^Be=wed^A=0$, so there exists a unique protomorphism $d^C\dd C\to C$ of
degree $-1$ such that $wd^B=d^Cw$. Also $d^Cd^Cw=d^Cwd^B=wd^Bd^B=0$. So $d^Cd^C=0$. Therefore $C$
becomes an object of $\mathrm{DG}\CA$ and $w$ a chain map of degree $0$. So $w$ is the cokernel of $x$
in $\mathrm{DG}\CA$. By Theorem~\ref{DGcase}, $\mathrm{DG}\CA$ is Cauchy complete.  
\end{proof}

\begin{center}
--------------------------------------------------------
\end{center}

\appendix

\end{document}